\documentclass[12pt]{amsart}

\usepackage{a4,xy,amssymb,diagrams,times,array,color,graphicx}

\usepackage[metapost,truebbox]{mfpic}
\opengraphsfile{paperfile}

\xyoption{poly}
\xyoption{arc}
\xyoption{knot}
\xyoption{curve}
\xyoption{arrow}
\xyoption{all}

\makeindex
 
\newcommand{\N}{\mathbb{N}}
\newcommand{\C}{\mathbb{C}}
\newcommand{\Z}{\mathbb{Z}}

\newcommand{\wis}[1]{{\text{\em \usefont{OT1}{cmtt}{m}{n} #1}}}
\newcommand{\Oscr}{\mathcal{O}}
\newcommand{\FF}{\mathbb{F}}
\newcommand{\V}{\mathbb{V}}
\newcommand{\Q}{\mathbb{Q}}

\newcommand{\vtx}[1]{*+[o][F-]{\scriptscriptstyle #1}}

\newcommand{\PSL}{\operatorname{PSL}}

\newcolumntype{C}{>{$}c<{$}}

\newcommand{\aar}[2]{ \ar@{-}|(.25){{\bf #1}}|{\bullet}|(.75){{\bf  #2}}}

\title{(non)commutative F-un geometry}
\author{Lieven Le Bruyn} 
\address{Department of Mathematics, University of Antwerp \\ 
 Middelheimlaan 1, B-2020 Antwerp (Belgium) \\ {\tt lieven.lebruyn@ua.ac.be}}

\begin{document}
\sloppy
 
\def\ldb{\mathopen{\{\!\!\{}} \def\rdb{\mathclose{\}\!\!\}}}


 \begin{abstract}
 Stressing the role of dual coalgebras, we modify the definition of affine schemes over the 'field with one element'. This clarifies the appearance of Habiro-type rings in the commutative case, and, allows a natural noncommutative generalization, the study of representations of discrete groups and their profinite completions being our main motivation.
\end{abstract}

\maketitle

\section{Commutative F-un geometry}

In this section we will recall the definition of affine schemes over the mythical field $\mathbb{F}_1$ with one element, originally due to Christophe Soul\'e \cite{Soule} and refined later by Alain Connes and Katia Consani \cite{CC}. This approach is based on functors from abelian groups to sets satisfying a universal property with respect to an integral- and a complex affine scheme. 
We will modify this definition slightly by replacing these affine schemes by integral- resp. complex dual coalgebras. This amounts to restricting to \'etale local data of the affine schemes and has the additional advantage that the definition can be extended verbatim to the noncommutative world as we will outline in the next section. Another advantage of the coalgebra approach is that it inevitably leads to the introduction of the Habiro ring \cite{Habiro} in the easiest example, that of the multiplicative group. This might be compared to recent work by Yuri I. Manin \cite{Manin} and Matilde Marcolli \cite{Marcolli}.

\subsection{} For a commutative ring $\wis{k}$ we will denote with $\wis{k-calg}$, resp. $\wis{k-alg}$, the category of all commutative $\wis{k}$-algebras, resp. the category of all $\wis{k}$-algebras. and with morphisms all $\wis{k}$-algebra morphisms. For two objects $A,B$ in $\wis{k-alg}$ we will denote the set of all $\wis{k}$-algebra morphisms from $A$ to $B$ by $(A,B)_{\wis{k}}$. 

\subsection{} Grothendieck introduced the category $\wis{k-caff}$ of all affine schemes living over a commutative ring $\wis{k}$ to be the category dual to the category $\wis{k-calg}$ of all commutative $\wis{k}$-algebras, that is, $\wis{k-caff} = (\wis{k-calg})^o$. 
One way to realize this duality is to associate to a commutative $\wis{k}$-algebra $A$ a covariant functor, {\em the functor of points} $\wis{h}_A$, 
\[
\wis{h}_A~:~\wis{k-calg} \rTo \wis{sets} \qquad B \mapsto (A,B)_{\wis{k}} \]
Alternatively, one can associate to $A$ a more classical geometric object, {\em the affine scheme} $\wis{spec}(A)$. This consists of a topological space $spec(A)$, the set of all prime ideals of $A$ equipped with the Zariski topology, together with a sheaf of rings $\Oscr_A$ on it, called the structure sheaf of $A$. The ring $A$ is recovered as the ring of global sections. Whereas both approaches are equivalent, it should be clear that the functorial point of view lends itself more easily to generalizations.

\subsection{} F-un or $\mathbb{F}_1$, the field with one element, is a virtual object which might be thought of as a 'ring' living under $\Z$. $\mathbb{F}_1$-believers base their f-unny intuition on the following two mantras :

\begin{itemize}
\item{$\mathbb{F}_1$ forgets about additive data and retains only multiplicative data.}
\item{$\mathbb{F}_1$-objects only acquire flesh when extended to $\Z$ (or $\C$).}
\end{itemize}

As an example, an $\mathbb{F}_1$-vectorspace is merely a set $V$ as there is no addition of vectors and just one element to use for scalar multiplication. Hence, the dimension of $V$ equals the cardinality of $V$ as a set. Next one should specify the classical objects one obtains after 'extending' $V$ to the integers or to the complex numbers. The correct integral version of a vectorspace is a lattice, so one defines $V \otimes_{\mathbb{F}_1} \Z$ to be the free $\Z$-lattice $\Z V$ on $V$. Analogously, one defines the extension of $V$ to the complex numbers, $V \otimes_{\mathbb{F}_1} \C$ to be the complex vectorspace $\C V$ with basis the set $V$.

But then, linear maps between $\mathbb{F}_1$-vectorspaces will be just set-maps and invertible maps are bijections, whence the group $GL_n(\mathbb{F}_1)$ is the symmetric group $S_n$. For a group $G$, an $n$-dimensional representation over $\mathbb{F}_1$ will then be a groupmorphism $\rho : G \rTo S_n$, that is, a permutation representation of $G$. Irreducible $G$-representations over $\mathbb{F}_1$ are then transitive permutation representations, and so on.

\subsection{} \label{commutative} In analogy with the finite field case, one expects there to be a unique $n$-dimensional field extension of $\mathbb{F}_1$ which we will denote by $\mathbb{F}_{1^n}$. This has to be a set with $n$ elements allowing a multiplication, whence the proposal to take $\mathbb{F}_{1^n} = C_n$ the cyclic group of order $n$. Extending $\mathbb{F}_{1^n}$ to the integers or complex numbers we should obtain a commutative algebra of rank resp. dimension $n$. Christophe Soul\'e \cite{Soule} proposed to take the integral- and complex group-algebras
\[
\mathbb{F}_{1^n} \otimes_{\mathbb{F}_1} \Z \simeq \Z C_n \quad \text{and} \quad \mathbb{F}_{1^n} \otimes_{\mathbb{F}_1} \C \simeq \C G \]
More generally, he proposed to take as the category of all commutative $\mathbb{F}_1$-algebras the category of all 
 finite (!) abelian groups, that is, $\mathbb{F}_1-\wis{calg} = \wis{abelian}$. For any abelian group $G$ we then have to make sense of the extended algebras which we take again to be the group-algebras 
 \[
 G \otimes_{\mathbb{F}_1} \Z \simeq \Z G \quad \text{and} \quad G \otimes_{\mathbb{F}_1} \C \simeq \C G \]
Having a notion for commutative $\mathbb{F}_1$-algebras,  Soul\'e takes Grothendieck functor of points approach to define {\em affine $\mathbb{F}_1$-schemes}. This should be a covariant functor
\[
X~:~\wis{abelian} \rTo \wis{sets} \]
connecting nicely to the functor of points of an affine integral- and complex-scheme. More precisely, Soul\'e \cite{Soule} and later Connes and Consani \cite{CC} require the following data
\begin{itemize}
\item{a complex affine commutative algebra $A \in \C-\wis{calg}$}
\item{an integral algebra $B \in \Z-\wis{calg}$ such that
$B \otimes_{\Z} \C \rInto A$}
\item{
a natural transformation $ev : X \rTo h_A$, called the 'evaluation' map}
\item{
an inclusion of functors $i : X \rInto h_B$}
\end{itemize}
satisfying the following universal property : given any integral algebra $C \in \Z-\wis{calg}$, any natural transformation $f : X \rTo h_C$ and any natural transformation $g : h_A \rTo h_{C \otimes_Z \C}$ making the upper square commute 
\[
\xymatrix{X \ar[rr]^{ev} \ar[rd]^f \ar[dd]^i & & h_A \ar[rd]^{g} & \\
& h_C \ar[rr]^{- \otimes \C} & & h_{C \otimes \C} \\
h_B \ar@{.>}[ru]^{\exists} \ar[rr]^{- \otimes \C} & & h_{B \otimes \C} \ar@{.>}[ru] &}
\]
there ought to be a natural transformation $h_B \rTo h_C$ making the entire diagram commute. This means that $\wis{spec}(B)$ is the best affine integral scheme approximating the functor $X$. Note that by Yoneda's lemma this means that one can reconstruct from the $\C$-algebra morphism $\psi : C \otimes \C \rTo A$ determining the natural transformation $g=- \circ \psi$ a $\Z$-algebra morphism $\phi : C \rTo B$ compatible with the inclusion $B \otimes \C \rInto A$. This means that for every abelian group $G$ we have a commuting diagram
\[
\xymatrix{X(G) \ar[rr]^{ev} \ar[rd]^f \ar[dd]^i & & (A,\C G)_{\C} \ar[rd]^{- \circ \psi} & \\
& (C,\Z G)_{\Z} \ar[rr]^{- \otimes \C} & & (C \otimes \C,\C G)_{\C} \\
(B,\Z G)_{\Z} \ar[ru]^{- \circ \phi} \ar[rr]^{- \otimes \C} & & (B \otimes \C,\C G)_{\C} \ar[ru] &}
\]

\subsection{} \label{example} The archetypical example being the multiplicative group. Consider the forgetful functor
\[
\mathbb{G}_m~:~\wis{abelian} \rTo \wis{sets} \qquad G \mapsto G \]
Take $A = \C[q^{\pm}]$ and $B = \Z[q^{\pm}]$, then their functors of points are exactly the multiplicative group scheme, that is give the groups of units
\[
h_A(D) = D^* \quad \text{and} \quad h_B(C)=C^* \]
for all $D \in \C-\wis{calg}$ and $C \in \Z-\wis{calg}$. We can then take both $i$ and $ev$ the natural transformation taking $F(G)=G$ to the subgroup of units $G \subset (\Z G)^* \subset (\C G)^*$. 

Remains only to prove the universal property. Let the natural transformation $g : h_{\C[q^{\pm}]} \rTo h_{C \otimes \C}$ be determined by the $\C$-algebra morphism $\psi : C \otimes \C \rTo \C[q^{\pm}]$ and let $N$ be a natural number larger than the degree of all $\psi(c)$ where $c$ is one of the $\Z$-algebra generators of $C$. Consider the finite cyclic group $C_N = \langle g \rangle$, then tracing the element $g$ around the above diagram gives the commutative diagram
\[
\begin{diagram}
C & \rInto & C \otimes \C & \rTo^{\psi} & \C[q^{\pm}] \\
\dTo^{\phi} & & & & \dOnto^{\pi} \\
\Z C_N & & \rInto & & \C C_N = \frac{\C[q^{\pm}]}{(q^N-1)} \end{diagram} \]
where $\phi = f(g)$. Repeating this argument,  $\pi(\psi(c)) = \psi(c) = \phi(c)$ for all $\Z$-generators of $C$, whence we have that $\psi(C) \subset \Z[q^{\pm}]$ giving the required natural transformation $h_{\Z[q^{\pm}]} \rTo h_C$.

\subsection{}  Observe that Soul\'e uses only finite abelian groups and hence we do not require the full functor of points, but rather the restricted functors
\[
\wis{h}'_A~:~\wis{k-fd.calg} \rTo \wis{sets} \qquad B \mapsto (A,B)_{\wis{k}} \]
where $\wis{k-fd.calg}$ is the category of all {\em finite dimensional} commutative $\wis{k}$-algebras. On the 'geometric' level we might still use the affine scheme $\wis{spec}(A)$ as this object contains more information than $\wis{h}'_A$, but we'd rather use a slimmer geometric object having the same amount of information as the restricted functor of points. It will turn out that the object we propose can be extended verbatim to the noncommutative world, whereas trying to extend affine schemes is known to lead to major difficulties.

\subsection{} Let us consider the complex case first. For $A \in \C-\wis{calg}$, we define the (finite) dual coalgebra $A^o$ to be the collection of all $\C$-linear maps $\lambda : A \rTo \C$ whose kernel contains a cofinite ideal $I \triangleleft A$. The dual maps to the multiplication and unit map of $A$ then define a coalgebra structure on $A^o$, see for example Sweedler's monograph \cite{Sweedler}. For $B$ a finite dimensional $\C$-algebra, any $\C$-algebra morphism $A \rTo B$ dualizes to a $\C$-coalgebra map $B^* \rTo A^o$ and as a coalgebra is the limit of its finite dimensional sub-coalgebras we see that the dual coalgebra $A^o$ contains the same information as the restricted functor of points $\wis{h}'_A$. We will now turn $A^o$ into our desired 'geometric' object.

As $A$ is commutative, any finite dimensional quotient $A/I \simeq L_{\mathfrak{m}_1} \oplus \hdots \oplus L_{\mathfrak{m}_k}$ splits into a direct sum of locals and hence the dual subcoalgebra $(A/I)^*$ is the direct sum of pointed coalgebras $(L_{\mathfrak{m}})^*$ which are subcoalgebras of the enveloping algebra of the abelian Lie-algebra of tangent-vectors $(\mathfrak{m}/\mathfrak{m}^2)^*$. Taking limits we have that
\[
A^o = \bigoplus_{\mathfrak{m} \in \wis{max}(A)} P_{\mathfrak{m}} \]
with $P_{\mathfrak{m}} \subset U((\mathfrak{m}/\mathfrak{m}^2)^*)$. In particular, we obtain the maximal ideals $\wis{max}(A)$ as the group-like elements of $A^o$, or equivalently, as the direct factors of the coradical $corad(A^o)$. Elements of $A$ naturally evaluate on $A^o$ (and hence on the coradical) and induce the usual Zariski topology on $\wis{max}(A)$.

We thus recover from the dual coalgebra $A^o$ the maximal ideal spectrum of $A$. But, $A^o$ contains a lot more local information. This is best seen by taking the full dual algebra $A^{o*}$ of $A^o$ giving rise to a Taylor-embedding (sending a function to its Taylor series expansions in all points)
\[
A \rInto A^{o*} = \prod_{\mathfrak{m} \in \wis{max}(A)} \hat{\Oscr}_{A,\mathfrak{m}} \]
where $\hat{\Oscr}_{A,\mathfrak{m}}$ is the $\mathfrak{m}$-adic completion of $A$ (that is the stalk of the structure sheaf in the \'etale topology). 

Concluding, the restricted functor of points $\wis{h}'_A$, or equivalently the dual coalgebra $A^o$, contains enough information to recover the analytic (or \'etale) local information in all the closed points of $\wis{spec}(A)$.

\subsection{} An affine F-un scheme $X : \wis{abelian} \rTo \wis{sets}$ connects to the complex picture via the evaluation natural transformation $ev : X \rTo \wis{h}'_A$. The discussion above leads to the introduction of an analytic ring of functions $\mathbb{F}_1[X]^{an}$ of which we now have a complex interpretation
\[
\mathbb{F}_1[X]^{an} \otimes_{\mathbb{F}_1} \C = \bigcap_{\mathfrak{m} \in Im(ev)} \hat{\Oscr}_{A,\mathfrak{m}} \]
With $Im(ev)$ we denote the images of all maps $\wis{max}(\C G) \rTo \wis{max}(A)$ coming from the algebra maps $A \rTo \C G$ contained in $ev(F(G)) \subset \wis{h}'_A(\C G)$.

For the example \ref{example} of the forgetful functor, we have $A = \C[q^{\pm}]$ and hence $\wis{max}(A) = \C^*$ and
\[
\C[q^{\pm}]^{o*} = \prod_{\alpha \in \C^*} \C[[q-\alpha]] \]
For any finite abelian group $G$, $\wis{max}(\C G)$ is the set of characters of $G$ and under the evaluation map an element $g \in F(G)=G$ maps a  character $\chi$ to its value $\chi(g)$, which are of course all roots of unity. Hence, if we vary over all finite abelian groups we obtain
\[
\mathbb{F}_1[q^{\pm}]^{an} \otimes_{\mathbb{F}_1} \C = \bigcap_{\lambda \in \mu_{\infty}} \C[[q-\lambda]] \]
Observe that $\mu_{\infty}$, the set of all roots of unity, is a Zariski dense set in $\wis{max}(\C[q^{\pm}]) = \C^*$.

\subsection{} Whereas the new complex picture based on the dual coalgebra is still pretty close to the usual affine scheme, this changes drastically  in the integral picture. For a $\Z$-algebra $B$ we have to consider the restricted functor of points
\[
\wis{h}'_B~:~\Z-\wis{fp.calg} \rTo \wis{sets} \qquad C \mapsto (B,C)_{\Z} \]
where $\Z-\wis{fp.calg}$ is the category of all commutative $\Z$-algebras which are finite projective $\Z$-modules. Again, this restricted functor contains the same information as the dual $\Z$-coalgebra
\[
B^o = \underset{\rightarrow}{lim}~Hom_{\Z}(B/I,\Z) \]
where the limit is taken over all ideals $I \triangleleft B$ such that $B/I$ is a projective $\Z$-module of finite rank. If we try to mimic the complex description of the dual coalgebra we are led to consider a certain subset of all coheight one prime ideals of $B$
\[
\wis{submax}(B) = \{ P \in spec(B)~|~\text{$B/P$ is a free $\Z$-module of finite rank} \} \]
Note that closed points in $\wis{spec}(B)$ are {\em not} contained in $\wis{submax}(B)$. Therefore we face the problem that different elements $P,P' \in \wis{submax}(B)$ are usually not comaximal and hence that we no longer have a direct sum decomposition of $B^o$ over this set (as was the case for the complex dual coalgebra).

As we will recall in the next section, we are familiar with such situations in noncommutative algebra, where even maximal ideals can belong to the same 'clique', that is, that the corresponding simple representations have nontrivial extensions. Using this noncommutative intuition, we therefore impose a clique-relation on the elements of $\wis{submax}(B)$
\[
P  \leftrightarrow P' \qquad \text{iff} \qquad P+P' \not= B \]
This relation should be thought of as a 'nearness' condition. Observe that any $P \in \wis{submax}(B)$ determines a finite collection of points in $\wis{max}(B \otimes_{\Z} \C)$ and hence we can extend this nearness relation on the points of $\wis{max}(B)$. Observe that this relation is clearly invariant under the action of the absolute Galois group $Gal(\overline{\Q}/\Q)$.

The different cliques determine the direct sum decomposition of the $\Z$-coalgebra $B^o$ and hence also of the Taylor-like ring of functions $B^{o*}$. Fully describing the dual $\Z$-coalgebra $B^o$ usually is a very difficult task and therefore, as in the complex case, when we are studying F-un geometry we restrict to that part determined by the elements in $Im(i)$ where $i~:~F \rTo \wis{h}'_B$ is the inclusion of functors determined by the affine F-un scheme $F~:~\wis{abelian} \rTo \wis{sets}$.

\subsection{} Let us consider again the example of the multiplicative group and indicate how the $\Z$-coalgebra approach leads to the introduction of the Habiro ring.

The ideals $I \triangleleft B= \Z[q^{\pm}]$ such that $B/I$ is a free $\Z$-module of finite rank are precisely the principal ideals $I = (f(q))$ where $f(q)$ is a monic polynomial. Hence,
\[
\wis{submax}(\Z[q^{\pm}]) = \{ (p(q))~:~\text{$p(q)$ is monic and irreducible} \} \]
Because $\Z[q^{\pm}]$ is a unique factorization domain we can decompose any monic polynomial uniquely into irreducible factors
\[
f(q) = p_1(q)^{n_1} \hdots p_k(q)^{n_k} \]
and we would like to use this fact, as in the complex case, to decompose the (linear duals) finite rank $\Z$-algebra quotients over $\wis{submax}(\Z[q^{\pm}])$. However, 
\[
\frac{\Z[q^{\pm}]}{(f(q))} \not= \frac{\Z[q^{\pm}]}{(p_1(q))^{n_1}} \oplus \hdots \oplus  \frac{\Z[q^{\pm}]}{(p_k(q))^{n_k}} \]
as the different primes $(p_i(q))$ and $(p_j(q))$ do not have to be comaximal. 
This problem makes it impossible to split the description of the dual coalgebra over the 'points' as in the complex case. Hence, we have no other option but to describe it as a direct limit
\[
\Z[q^{\pm}]^o = \underset{\rightarrow}{lim}~(\frac{\Z[q^{\pm}]}{(f(q))})^* \]
where the limit is considered with respect to divisibility of polynomials as there are natural inclusions of $\Z$-coalgebras
\[
(\frac{\Z[q^{\pm}]}{(f(q))})^* \rInto (\frac{\Z[q^{\pm}]}{(g(q))})^* \qquad \text{whenever} \qquad f(q) | g(q) \]
As in the complex case we are then interested in the dual algebra of $\Z[q^{\pm}]^o$ and the natural algebra map
\[
\Z[q^{\pm}]^o \rInto (\Z[q^{\pm}]^o)^* = \underset{\leftarrow}{lim}~\frac{\Z[q^{\pm}]}{(f(q))} \]
and it is clear that in the description of the algebra on the right-hand side completions at principal ideals will constitute a main ingredient.

While we can do all these calculations to some extend, we are primarily interested in that part of $\wis{submax}(\Z[q^{\pm}])$ in the image of the inclusion functor, that is 
\[
Im(i) = \N = \{ (\Phi_1(q)), (\Phi_2(q)), \hdots , (\Phi_n(q)), \hdots \} \subset  \wis{submax}(\Z[q^{\pm}]) \]
We will confuse the natural number $n$ with the corresponding cyclotomic polynomial $\Phi_n(q)$ or with the height one prime generated by it. With this identification $\N$ is the integral analog of the set of all roots of unity $\boldsymbol{\mu}_{\infty}$ in the complex case.

In the case of cyclotomic polynomials we have complete information about possible co-maximality \begin{itemize}
\item{If $\frac{m}{n} \not= p^k$ for some prime number $p$, then $(\Phi_m(q),\Phi_n(q))=1$ that is the cyclotomic prime ideals are comaximal.}
\item{If $\frac{m}{n}=p^k$ for some prime number $p$, then $\Phi_m(q) \equiv \Phi_n(q)^d~\wis{mod}~(p)$ for some integer $d$, hence the cyclotomic primes are not comaximal.}
\end{itemize}
Therefore, the relevant clique-relation is
\[
n \leftrightarrow m \qquad \text{if and only if} \qquad \frac{m}{n}=p^{\pm k} \]
inducing on the complex level the $Gal(\overline{\Q}/\Q)$-invariant nearness condition on roots of unity $\lambda, \mu \in \mu_{\infty}$
\[
\lambda \leftrightarrow \mu \qquad \text{iff} \qquad \frac{\lambda}{\mu}~\text{is of order $p^k$} \]
for some prime number $p$. 

Yuri I. Manin argues in \cite{Manin} that we should take the analogy between the integral affine scheme $\wis{spec}(\Z[q^{\pm}])$ and the (complex) affine plane more seriously and that, besides the arithmetic axis, one should also consider a projection to the 'geometric axis' (which should then be viewed as the affine $\mathbb{F}_1$-scheme corresponding to $\mathbb{F}_1[q^{\pm}]$. He proposed that the zero sets of the cyclotomic polynomials $\Phi_n(q)$ for all integers $n$ should be considered as the union of the fibers in this second projection. That is, we should have the following picture :

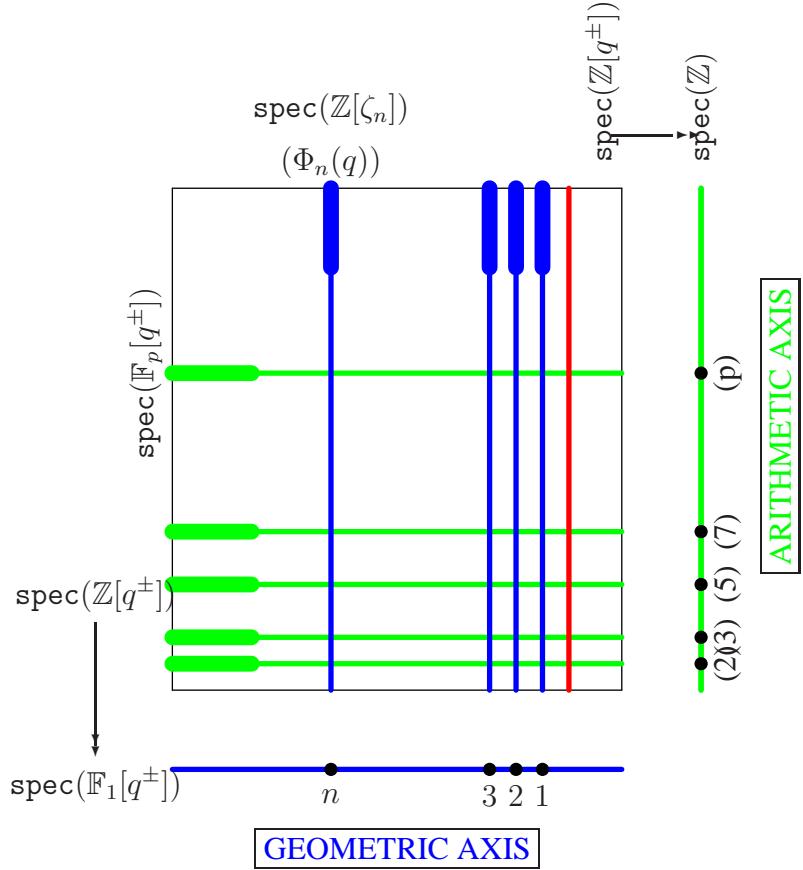
\[
\rotatebox{90}{
\begin{mfpic}[1]{0}{0}{100}{100}
\lines{(10,30),(200,30)}
\lines{(10,30),(10,200)}
\lines{(10,200),(200,200)}
\lines{(200,30),(200,200)}
\pen{2pt}
\drawcolor{green}
\lines{(10,0),(200,0)}
\lines{(20,30),(20,200)}
\lines{(30,30),(30,200)}
\lines{(50,30),(50,200)}
\lines{(70,30),(70,200)}
\lines{(130,30),(130,200)}
\lines{(10,0),(200,0)}
\pen{6pt}
\lines{(20,170),(20,200)}
\lines{(30,170),(30,200)}
\lines{(50,170),(50,200)}
\lines{(70,170),(70,200)}
\lines{(130,170),(130,200)}
\pen{2pt}
\drawcolor{red}
\lines{(10,50),(200,50)}
\drawcolor{blue}
\lines{(10,70),(200,70)}
\lines{(10,60),(200,60)}
\lines{(10,80),(200,80)}
\lines{(10,140),(200,140)}
\lines{(-20,30),(-20,200)}
\pen{6pt}
\lines{(170,70),(200,70)}
\lines{(170,60),(200,60)}
\lines{(170,80),(200,80)}
\lines{(170,140),(200,140)}
\pen{2pt}
\pointsize=5pt
\point{(20,0),(30,0),(50,0),(70,0),(130,0),(-20,60),(-20,70),(-20,80),(-20,140)}
\tlabel[cc](20,-10){(2)}
\tlabel[cc](30,-10){(3)}
\tlabel[cc](50,-10){(5)}
\tlabel[cc](70,-10){(7)}
\tlabel[cc](130,-10){(p)}
\tlabel[cc](130,210){$\wis{spec}(\mathbb{F}_p[q^{\pm}])$}
\tlabel[cl](210,37){$\wis{spec}(\Z[q^{\pm}])$}
\tlabel[cc](220,20){$\begin{diagram} \\ \dOnto  \end{diagram}$}
\tlabel[cl](210,0){$\wis{spec}(\Z)$}
\tlabel[cc](110,-30){\fbox{\textcolor{green}{ARITHMETIC AXIS}}}
\tlabel[cc](-50,115){\rotatebox{270}{\fbox{\textcolor{blue}{GEOMETRIC AXIS}}}}
\tlabel[cc](-30,60){\rotatebox{270}{$1$}}
\tlabel[cc](-30,70){\rotatebox{270}{$2$}}
\tlabel[cc](-30,80){\rotatebox{270}{$3$}}
\tlabel[cc](-30,140){\rotatebox{270}{$n$}}
\tlabel[cc](210,140){\rotatebox{270}{$(\Phi_n(q))$}}
\tlabel[cc](230,140){\rotatebox{270}{$\wis{spec}(\Z[\zeta_n])$}}
\tlabel[cc](10,230){\rotatebox{270}{$\begin{diagram} \wis{spec}(\Z[q^{\pm}]) \\ \dOnto \\ \wis{spec}(\mathbb{F}_1[q^{\pm}]) \end{diagram}$}}
\end{mfpic}}
\]

Note that this is an over-simplification. Whereas the different green fibers for the projection to the arithmetic axis are clearly comaximal, the blue fibers are not. For example, the zero sets $\V(\Phi_2(q))$ and $\V(\Phi_1(q))$ share the maximal ideal $(2,q-1)$. The clique-relation encodes how the blue fibers intersect each other.

The clique-relation is important to relate different completions occurring in the F-un determined part of the algebra $(\Z[q^{\pm}]^o)^*$ as was proved by Kazuo Habiro \cite{Habiro}. Let us define for any subset $S \subset \N$ the completion
\[
\Z[q^{\pm}]^S = \underset{\underset{p \in \Phi_S^*}{\leftarrow}}{lim}~\frac{\Z[q^{\pm}]}{(p)} \]
where $\Phi_S^*$ is the set of monic polynomials generated by all $\Phi_n(q)$ for $n \in S$. Among the many precise results proved in \cite{Habiro} we mention these two

\begin{enumerate}
\item{If $S' \subset S$ and if every clique-component of $S$ contains an element from $S'$, then the natural map is an inclusion
\[
\rho^S_{S'}~:~\Z[q^{\pm}]^S \rInto \Z[q^{\pm}]^{S'} \]}
\item{If $S$ is a saturated subset of $\N$ meaning that for every $n \in S$ also its divisor-set $\langle n \rangle = \{ m | n \}$ is contained in $S$, then
\[
\Z[q^{\pm}]^S = \bigcap_{n \in S} \Z[q^{\pm}]^{\langle n \rangle} = \bigcap_{n \in S} \widehat{\Z[q^{\pm}]}_{(q^n-1)} \]
where the terms on the right-hand side are the $I$-adic completions where $I=(q^n-1)$.}
\end{enumerate}
Using these properties it is then natural to define the integral version of the ring of analytic functions on the multiplicative group scheme over $\FF_1$ to be
\[
\FF_1[q^{\pm}]^{an} \otimes_{\FF_1} \Z \simeq \bigcap_{n \in \N} \widehat{\Z[q^{\pm}]}_{(q^n-1)} = \Z[q^{\pm}]^{\N} \]
This ring has a description very similar to that of the profinite integers replacing factorials by q-factorials
\[
\Z[q^{\pm}]^{\N} = \underset{\underset{n}{\leftarrow}}{lim}~\frac{\Z[q^{\pm}]}{((q^n-1)(q^{n-1}-1) \hdots (q-1))} \]
and as such its elements have a unique description as formal Laurent polynomials over $\Z$ of the form
\[
\sum_{n=0}^{\infty} a_n(q) (q^n-1)(q^{n-1}-1)\hdots (q-1) \in \Z[[q^{\pm}]] \qquad \text{with} \qquad deg(a_n(q)) < n \]
We observe that any such formal power series can be evaluated at a root of unity. Some elements of $\Z[q^{\pm}]^{\N}$ have been discovered before. For example, Maxim Kontsevich observed in his investigations on Feynman integrals that the formal power series
\[
\sum_{n=0}^{\infty} (1-q)(1-q^2) \hdots (1-q^n) \]
has a properly defined value in every root of unity. Subsequently, Don Zagier \cite{Zagier} proved the strange equality
\[
\sum_{n=0}^{\infty} (1-q)(1-q^2) \hdots (1-q^n) = - \frac{1}{2} \sum_{n=1}^{\infty} n \chi(n) q^{(n^2-1)/24} \]
where $\chi$ is the quadratic character of conductor $12$. The strange fact about this equality is that the two sides never make sense simultaneously. The left hand side diverges for all points within the unit circle and outside the unit circle and can be evaluated at roots of unity whereas the right hand side converges only within the unit circle and diverges everywhere else. What Zagier meant by this equality is that for all $\alpha \in \boldsymbol{\mu}_{\infty}$ the evaluation of the left hand side coincides with the radial limit of the function on the right hand side. Don Zagier says that the function on the right 'leak through roots of unity'.

\section{Noncommutative F-un geometry}

In this section we will extend Soul\'e's definition of an affine $\mathbb{F}_1$-scheme to the noncommutative case. Our main motivation is the study of finite dimensional representations of discrete groups, such as the braid groups or the modular group. We have seen that irreducible finite dimensional $\mathbb{F}_1$-representations of a group $\Gamma$ are exactly the finite transitive permutation representations $\Gamma/\Lambda$ where $\Lambda$ is of finite index in $\Gamma$. That is, all finite dimensional $\mathbb{F}_1$-representation theory of $\Gamma$ comes from its profinite completion $\hat{\Gamma} = \underset{\leftarrow}{lim}~\Gamma/\Lambda$, the limit taken over all finite index normal subgroups.

In the previous section we have worked out the special case when $\Gamma = \Z$. Here, the simple representations of $\hat{\Z}$ are the roots of unity $\mu_{\infty}$ and they are Zariski closed in all simples $\C^* = \wis{simp}(\Z)$. The clique-relation on $\mu_{\infty}$ was compatible with the action of the absolute Galois group and the Habiro ring 'feels' the inclusion $\mu_{\infty} \subset \C^*$, that is it contains the tangent information in a Galois-compatible way.

Here we extend some of these results to the case of a non-Abelian discrete group $\Gamma$ satisfying the property $\bullet$ : for every finite collection of elements $\{ g_1,\hdots,g_k \} \subset \Gamma$ there is a finite index subgroup $\Lambda \subset \Gamma$ such that the natural projection map gives an embedding $\{ g_1,\hdots,g_k \} \rInto \Gamma/\Lambda$. We will prove that such groups determine a noncommutative affine $\mathbb{F}_1$-scheme, the F-un information being given by the finite dimensional permutation representations, or equivalently, the representation theory of the profinite completion $\hat{\Gamma}$. We will show that $\wis{simp}(\hat{\Gamma})$ is Zariski dense in $\wis{simp}(\Gamma)$ and compute the tangent information of this embedding. That is, to a finite dimensional permutation representation $P=\Gamma/\Lambda$ we will associate a noncommutative gadget (a quiver, relations and a dimension vector) encoding all possible deformations of $P$ which are still $\Gamma$-representations. In relevant situations, including the case when $\Gamma$ is the modular group $\PSL_2(\Z)$ (in which case the permutation representations are Grothendieck's 'dessins d'enfants') some subsidiary noncommutative gadgets can be derived from this tangent information, such as the necklace Lie algebra \cite{LBBocklandt} and the singularity type \cite{LBBocklandtSymens}. It is to be expected that most of these noncommutative gadgets associated to dessins are in fact Galois invariants.

\subsection{} If we take commutative $\mathbb{F}_1$-algebras to be abelian groups, it make sense to identify the category of all $\mathbb{F}_1$-algebras with $\wis{groups}$ the category of all finite groups. Likewise, we have to extend Grothendieck's functor of points to all, that is including  also noncommutative, algebras. With these modifications we can extend Soul\'e's definition to the noncommutative world.

Define an affine noncommutative $\mathbb{F}_1$-scheme to be a covariant functor
\[
X~:~\wis{groups} \rTo \wis{sets} \]
from the category $\wis{groups}$ of all finite groups to $\wis{sets}$. We require that there is an affine $\C$-algebra $A$ and an evaluation natural transformation $ev : X \rTo \wis{h}_A = (A,-)_{\C}$, giving for every finite group $G$ an evaluation map $X(G) \rTo (A,\C G)_{\C}$. 
Moreover, there should be a 'best' integral affine algebra $B$ with an inclusion of functors $X \rInto \wis{h}_B = (B , -)_{\Z}$. 

That is, for every finite group $G$ we have an inclusion $X(G) \rInto (B,\Z G)_{\Z}$. Here, 'best' means that for every $\Z$-algebra $C$ and every natural transformation $X \rTo \wis{h}_C = (C,-)_{\Z}$ and every $\C$-algebra morphism $\psi : \C \otimes C \rTo A$ making the upper square in the diagram below commute for every finite group $G$
\[
\xymatrix{X(G) \ar[rr]^{ev} \ar[rd]^f \ar[dd]^i & & (A,\C G)_{\C} \ar[rd]^{g=- \circ \psi} & \\
& (C,\Z G)_{\Z} \ar[rr]^{- \otimes \C} & & (C \otimes \C, \C G)_{\C} \\
(B,\Z G)_{\Z} \ar@{.>}[ru]^{\exists - \circ \phi} \ar[rr]^{- \otimes \C} & & (B \otimes \C,\C G)_{\C} \ar@{.>}[ru] &}
\]
there exists a $\Z$-algebra morphism $\phi : C \rTo B$ making the entire diagram commute.

\subsection{} \label{dessins} Our first example of a noncommutative F-un scheme is Grothendieck's theory of 'dessins d'enfants'. Let $X_{\C}$ be a Riemann surface (projective algebraic curve) defined over $\overline{\Q}$, then Belyi proved that there is a degree $d$ map $\pi : C \rOnto \mathbb{P}^1_{\C}$ ramified only in the points $\{ 0,1,\infty \}$. The open interval $] 0,1 [$ lifts to $d$ intervals on $C$. The endpoints of different lifts can be identified on $X$ indicating how the different sheets should be glued together in a neighborhood of the ramification point. The resulting graph with $d$ edges on $C$ is then called the {\em dessin} of $C$ and as the absolute Galois group $Gal(\overline{\Q}/\Q)$ acts on the collection of all such curves, it also acts on the dessins. Writing out this action allows one to gain insight in the absolute Galois group. Hence it is a very important problem to find new Galois invariants of dessins.

We will be particularly interested in {\em modular} dessins, that is such that the preimages of $0$ all have valency 1 or 2 and the preimages of 1 all have valency 1 or 3 in the graph. Alternatively, this means that the curve can be viewed as the compactification of a quotient $C = \mathbb{H}/\Lambda$ of the upper-halfplane under the action of a subgroup $\Lambda$ of finite index in the modular group $\Gamma = PSL_2(\Z)$. That is, modular dessins are equivalent to finite dimensional permutation representations of the modular group. Therefore, one is interested in the functor
\[
X~:~\wis{groups} \rTo \wis{sets} \qquad G \mapsto G_{(2)} \times G_{(3)} \]
sending a group to the set of all permutation representations of $\Gamma$ determined by elements of $G$. As $\Gamma \simeq C_2 \ast C_3$ is the free product of a cyclic group of order 2 with a cyclic group of order 3, this functor sends a finite group $G$ to the set product of its elements of order 2 with the elements of order 3 : $G_{(2)} \times G_{(3)}$. This functor determines a noncommutative affine $\mathbb{F}_1$-scheme as we can take as the complex- and integral group-algebras
\[
A = \C \Gamma \quad \text{and} \quad B = \Z \Gamma \]
of the modular group. As any $\C$-algebra morphism $A = \C \Gamma \rTo \C G$ is determined by the images of the order two (resp. three) generators $x$ and $y$ we can take as the evaluation and inclusion maps
\[
ev~:~G_{(2)} \times G_{(3)} \rTo (\C \Gamma, \C G)_{\C} \qquad (g_2,g_3) \mapsto \begin{cases} x \mapsto g_2 \\ y \mapsto g_3 \end{cases} \]
\[
i~:~G_{(2)} \times G_{(3)} \rInto (\Z \Gamma, \Z G)_{\Z} \qquad (g_2,g_3) \mapsto \begin{cases} x \mapsto g_2 \\ y \mapsto g_3 \end{cases} \]
We can repeat the argument of~\ref{example} verbatim to prove that these data indeed define a noncommutative $\mathbb{F}_1$-scheme using the fact that the modular group $\Gamma$ satisfies condition $\bullet$.

\subsection{} \label{TQFT} The second example is motivated by 2-dimensional TQFT. To a Riemann surface $C$ of genus $g$ and any finite group $G$ one associates as topological invariant $Z_G(C)$ the number of fields on $C$ with gauge group $G$, or equivalently, the number of $G$-covers on $C$. By Frobenius-Schur this number is equal to
\[
Z_G(C) = \sum_{\chi} (\frac{| G |}{dim~\chi})^{2g-2} \]
where the sum runs over all irreducible representations $\chi$ of the finite group $G$. As the number of $G$-covers is equal to the number of group-morphisms $\pi_1(C) \rTo G$ from the fundamental group $\pi_1(C) = \langle x_1,\hdots,x_g,y_1,\hdots,y_g \rangle / (\prod x_iy_ix_i^{-1} y_i^{-1} )$, this motivates the functor
\[
X~:~\wis{groups} \rTo \wis{sets} \quad G \mapsto \{ (a_1,\hdots,a_g,b_1,\hdots,b_g) \in G^{2g}~:~\prod a_ib_ia_i^{-1}b_i^{-1} = 1 \} \]
This functor is again an affine noncommutative $\mathbb{F}_1$-scheme as we can take the integral- and complex group-algebras $A = \C \pi_1(C)$ and $B= \Z \pi_1(C)$ and the natural evaluation and inclusion maps. Once again, the defining "bestness" property is verified using the fact that $\pi_1(C)$ satisfies condition $\bullet$.

Also in this example, the $\mathbb{F}_1$-info is given by all finite permutation representations of the fundamental group $\pi_1(C)$. That is, the F-un information is contained in the profinite completion $\widehat{\pi_1(C)}$.

\subsection{}  These two examples illustrate that any discrete group $\Gamma$ satisfying condition $\bullet$ determines a noncommutative affine $\mathbb{F}_1$-scheme. The corresponding functor assigns to a group $G$ the set of all groupmorphisms $\Gamma \rTo G$ and takes as the complex- and integral algebras the complex and integral group-algebra of $\Gamma$.

As in the commutative case we do not require the full strength of the functor of points $\wis{h}_A : \wis{k-alg} \rTo \wis{sets}$ for a given (not necessarily commutative) $\wis{k}$-algebra $A$, but it suffices, for applications to F-un geometry, to restrict to finite dimensional $\wis{k}$-algebras
\[
\wis{h}'_A~:~\wis{k-fd.alg} \rTo \wis{sets} \qquad C \mapsto (A,C)_{\wis{k}} \]
If $\wis{k}$ is a field, the information contained in this restricted functor of points is equivalent to that contained in the dual coalgebra $A^o$. For this reason we want to associate noncommutative geometric data (say, a topological space and function) to the dual $\C$-coalgebra $A^o$ where $A$ is the complex algebra determining the evaluation natural transformation $ev : X \rTo \wis{h}'_A$.

Observe that in \cite{LBdualcoalgebra} we initiated the description of the dual coalgebra of any affine $\C$-algebra $A$ in terms of the $A_{\infty}$-structure on the Yoneda space of all finite dimensional simple $A$-representations. For the applications we have in mind here, that is, virtually free groups $G$ (such as the modular group $\Gamma = PSL_2(\Z)$), for which the group algebras $\C G$ is formally smooth by \cite{LBqurves}, or 2-Calabi-Yau algebras such as $\C \pi_1(C)$, we do not require the full power of $A_{\infty}$-theory and can give, at least in principle,  an explicit description of the dual coalgebra.

The geometric space associated to an affine $\C$-algebra $A$ will be the set of isomorphism classes of finite dimensional $A$-representations, which as in the commutative case, is the set of direct summands of the coradical of the dual coalgebra
\[
\wis{simp}(A) = corad(A^o) \]
In \cite{LBdualcoalgebra} we introduced a Zariski topology on $\wis{simp}(A)$ in terms of the measuring $A^o \otimes A \rTo \C$. Here we will follow a slightly different approach based on noncommutative functions.

For a $\C$-algebra $A$ we define the {\em noncommutative functions} to be the $\C$-vectorspace quotients
\[
\wis{functions}(A) = \mathfrak{g}_A = \frac{A}{[A,A]_{vect}} \]
where $[A,A]_{vect}$ is the subvectorspace (and {\em not} the ideal) spanned by all commutators in $A$. 
Note that in the classical case where $A = \C[X]$ is the commutative coordinate ring of an affine variety $X$, there is nothing to divide out and hence in this case we recover the coordinate ring $\mathfrak{g}_A = \C[X]$. If $A = \C G$ the group-algebra of a finite group $G$, then $\mathfrak{g}_A$ is the space dual to the space of character-functions of $G$. Hence, in both cases the linear functionals $\mathfrak{g}^*$ suffice to separate the {\em points} of $A$, that is $\wis{simp}(A)$.  We will show that for a general affine $\C$-algebra $A$ we do indeed have an embedding
\[
 \wis{simp}(A) \rInto \mathfrak{g}^* \]
 Consider the (commutative) affine scheme $\wis{rep}_n A$ of all $n$-dimensional representations. A quick and dirty way to describe its coordinate ring $\C [\wis{rep}_n A]$ is to take a finite set of algebra generators $\{ a_1,\hdots,a_m \}$ of $A$, consider a set of $mn^2$ commuting variables $\{ x_{ij}(k) : 1 \leq i,j \leq n, 1 \leq k \leq m \}$ and consider the ideal $I_n(A)$ of the polynomial algebra $\C[x_{ij}(k)~:~i,j,k]$ generated by all entries of all $n \times n$ matrices $f(X_1,\hdots,X_m)$ where $f(a_1,\hdots,a_m)$ runs over all relations holding in $A$ and where $X_k$ is the generic $n \times n$ matrix $(x_{ij}(k))_{i,j}$. Then,
\[
\C[\wis{rep}_n A] = \frac{\C [ x_{ij}(k)~:~i,j,k ]}{I_n(A)} \]
On the affine scheme $\wis{rep}_n A$ there is a natural action of $GL_n$, the orbits of which correspond exactly to the isomorphism classes of $n$-dimensional $A$-representations. Basic GIT-stuff tells us that one can classify the {\em closed} orbits by points of the quotient-scheme $\wis{iss}_n A = \wis{rep}_n A / GL_n$ corresponding to the affine ring of invariants
\[
\C[\wis{iss}_n A] = \C[\wis{rep}_n A]^{GL_n} \]
and Artin proved that the closed orbits are precisely the isoclasses of {\em semi-simple} representations.
 
Let us bring in our quotient $\mathfrak{g} _A= \frac{A}{[A,A]_{vect}}$. We can evaluate its elements on all points of $\wis{rep}_n A$ by {\em taking traces}. That is, each $g \in \mathfrak{g}$ defines a function
\[
\wis{rep}_n A \rTo \C \qquad M \mapsto tr(g)(M) \]
That is, lift $g$ to an element $a \in A$, write $a=f(a_1,\hdots,a_m)$ in terms of its generators, then if $(m_1,\hdots,m_k)$ are the matrices describing the $n$-dimensional representation $M$, then we define
\[
tr(g)(M) = Tr(f(m_1,\hdots,m_k)) \]
where $Tr$ is the standard trace map on $M_n(\C)$. Observe that this does not depend on the chosen lift $a$ as all traces of elements from $[A,A]_{vect}$ vanish. Observe that via this trace-trick we can view elements of $\mathfrak{g}^*$ indeed as {\em generalized characters} as each representation defines a linear functional
\[
\chi_M~:~\mathfrak{g} \rTo \C \qquad g \mapsto tr(g)(M) \]

It is a classical result that the ring of invariants $\C[\wis{rep}_n A]^{GL_n}$ is generated by the invariant functions $tr(g)$ when $g$ runs over $\mathfrak{g}$. So, indeed, linear functionals on $\mathfrak{g}$ do separate $n$-dimensional semi-simple representations (whence a fortiori also simples). Actually, we only showed separation of simples for a fixed $n$, but clearly one recovers the dimension from $tr(1)$. 
That is, we have proved that for any affine $\C$-algebra $A$, the generalized character values give an embedding
\[
\wis{simp}(A) \rInto \mathfrak{g}^*_A \]
We will make the set $\wis{simp}(A)$ into a topological space by taking as the basic opens 
\[
\mathbb{X}(g,\lambda) = \{ S \in \wis{simp} A~|~\chi_S(g) \not= \lambda \} \]
for all $g \in \mathfrak{g}_A$ and all $\lambda \in \C$. For example, all simples of dimension $n$ form a closed subset. The obtained topology we will call the {\em Zariski topology} on $\wis{simp}(A)$.

Our use of this topology is to prove a denseness result similar to the fact that roots of unity $\mu_{\infty}$ are Zariski dense in $\C^*$. Let $G$ be a discrete group, as every finite dimensional $\hat{G}$ representation factors over a finite group quotient of $G$ (and hence is semi-simple) we deduce that the dual coalgebra $(\C G)^o$ is co-semi-simple and hence
\[
\wis{simp}(\C \hat{G}) = (\C G)^o = corad((\C G)^o) \]
We claim that when $G$ is a discrete group satisfying condition $\bullet$, then
\[
\overline{\wis{simp}(\C \hat{G})} = \wis{simp}(\C G) \]
That is, the subset of simple representations of the profinite completion is Zariski dense in the noncommutative space $\wis{simp}(\C G)$. Observe that in the two examples given before, $\wis{simp}(\C \hat{G})$ is the image of the evaluation map determined by the F-un geometry, hence this result is a direct generalization of the commutative situation for the multiplicative group.

To prove this claim observe  that the space of noncommutative functions $\mathfrak{g} = \mathfrak{g}_{\C G}$ has as $\C$-basis the conjugacy classes of elements of $G$. Hence, any linear functional $\chi \in \mathfrak{g}^*$ is a linear combination
\[
\chi = \lambda_1 \chi_1 + \hdots + \lambda_k \chi_k \]
where the $\chi_i$ are character functions corresponding to distinct conjugacy classes of $G$. Vanishing of $\chi$ on the whole of $\wis{simp}(\C \hat{G})$ would imply that the characters $\lambda_1,\hdots,\lambda_k$ are linearly dependent on every finite quotient $G/H$, which is impossible by the assumption on $G$.

\subsection{} Let us recall briefly the main result of \cite{LBdualcoalgebra} describing the dual coalgebra $A^o$ of a general affine $\C$-algebra $A$ and indicate the geometric information contained in it. Let $Q$ be a possibly infinite quiver and $\C Q$ the vectorspace spanned on all paths in $Q$ of positive length. Then $\C Q$ is given a coalgebra structure (the {\em path coalgebra})
\[
\Delta(p) = \sum_{p=p_1.p_2} p_1 \otimes p_2 \qquad \epsilon(p) = \delta_{p,vertex} \]
where $p_1.p_2$ is the concatenation of paths and the counit maps non-vertex paths to zero. 

Starting from $A$ we will construct a huge quiver $Q_A$ having as its vertices the isoclasses of finite dimensional simple representations and with the number of arrows between them
\[
\#(S \rTo  S' ) = dim_{\C}~Ext^1_A(S,S') \]
We will now describe a certain subcoalgebra of the path coalgebra $\C Q_A$ and as any coalgebra is the direct limit of its finite dimensional subcoalgebras we may restrict attention to a finite collection of simples and consider the semi-simple representation $M = S_1 \oplus \hdots \oplus S_k$ with restricted path-coalgebra $\C Q_A | M$. There is a natural $A_{\infty}$-algebra structure on the Yoneda Ext-algebra $Ext^{\bullet}_A(M,M)$, in particular there are higher multiplication maps
\[
m_i~:~\underbrace{Ext^1_A(M,M) \otimes \hdots \otimes Ext^1_A(M,M)}_i \rTo Ext^2_A(M,M) \]
defining a linear map, called the homotopy Maurer-Cartan map
\[
HMC_M = \oplus_i m_i~:~\C Q_A | M \rTo Ext^2_A(M,M) \]
The main result of \cite{LBdualcoalgebra} asserts that the dual coalgebra $A^o$ is Morita-Takeuchi equivalent to the largest subcoalgebra of $\C Q_A$ contained in the kernel of $HMC_M$ for all semi-simple representations $M$.

We will now describe the geometric content of the dual coalgebra. Recall that in the commutative case we had that the full linear dual of the dual coalgebra $(\C[X]^o)^* = \prod_x \hat{\Oscr}_{X,x}$ gave us back all the completed local rings at points of $X$. In the general case, assume as above that $M = S_1 \oplus \hdots \oplus S_k$ is a semi-simple representation with all simple factors distinct.The action of $A$ on $M$ gives rise to an epimorphism
\[
A \rOnto^{\pi_M} B_M=M_{n_1}(\C) \oplus \hdots \oplus M_{n_k}(\C) \]
and let us denote $\mathfrak{m} = Ker(\pi_M)$. If $C_M$ is the maximal subcoalgebra of $\C Q_A | M$ contained in the kernel of the $HMC_M$, then we can generalize the commutative situation as follows. The $\mathfrak{m}$-adic completion of $A$ is Morita equivalent to the full linear dual of $C_M$
\[
\hat{A}_{\mathfrak{m}} \sim_M (C_M)^* \]
This means that all $\mathfrak{m}$-adic completion of $A$ can be computed from the dual coalgebra $A^o$ and that each of them is a ring Morita equivalent to (the completion of) a path algebra of the quiver $(Q_A | M)^*$ modulo certain relations coming from the $A_{\infty}$-structure.

\subsection{} Recall that a $\C$-algebra $A$ is said to be {\em smooth} if and only if the kernel of the multiplication map
\[
\Omega^1_A = Ker( A \otimes A \rTo^m A) \]
is a projective $A$-bimodule. Because $Ext^2_A(M,N)=0$ for all finite dimensional $A$-representations when $A$ is smooth, we have from the above general result that the $\mathfrak{m}$-adic completion $\hat{A}_{\mathfrak{m}}$ is Morita-equivalent to the completion of the path algebra $\C (Q_A | M)^*$ where we recall that this quiver depends only on the dimensions of the ext-groups $Ext^1_A(S_i,S_j)$.

In fact, in this case we do not have to use the full strength of the general result and deduce this fact from the formal neighborhood theorem for smooth algebras due to Cuntz and Quillen \cite[\S 6]{CuntzQuillen}. Note that $Ker(\pi_M)= \mathfrak{m}$ has a natural $B=B_M$-bimodule structure. In analogy with the Zariski tangent space in the commutative case, we define
\[
T_M = \left( \frac{\mathfrak{m}}{\mathfrak{m}^2} \right)^{*} \]
Because $B$ is a semi-simple algebra  the simple $B$-bimodules are either of the form $M_{n_i}(\C)$ (with trivial action of the other components of $B$) or $M_{n_i \times n_j}(\C)$ with the component $M_{n_i}(\C)$ (resp. $M_{n_j}(\C)$) acting by left (resp. right) multiplication and all other actions being trivial. That is, there is a natural one-to-one correspondence between
\[
\wis{bimod}~B \leftrightarrow \wis{quiver}_n \]
isoclasses of $B$-bimodules and quivers $n$ vertices (the number of simple components). Under this correspondence, $B$-bimodule duals corresponds to taking the opposite quiver. Hence, the tangent space $T_M$ can be identified with a quiver on the vertices $\{ S_1,\hdots,S_n \}$ which we will now show is the opposite quiver of $Q_A | M$.

By the formal tubular neighborhood theorem of Cuntz and Quillen \cite[\S 6]{CuntzQuillen} (using the fact that semi-simple algebras are formally smooth) we have an isomorphism of completed algebras between the $\mathfrak{m}$-adic completion of $A$
\[
\hat{A}_{\mathfrak{m}} = \underset{\leftarrow}{lim}~A/\mathfrak{m}^n \]
where $\mathfrak{m} = Ker(\pi)$ as above, and, the completion (with respect to the natural gradation) of the tensor-algebra $T_B(\mathfrak{m}/\mathfrak{m}^2)$. That is, when we view $T_M$ as a quiver, then there is a Morita-equivalence
\[
\hat{A}_{\mathfrak{m}} \underset{M}{\sim} \widehat{\C T_M^{\vee}} \]
between the completion $\hat{A}_{\mathfrak{m}}$ and the completion (with respect to the gradation giving all arrows degree one) of the path-algebra $\C T_M^{\vee}$ of the opposite quiver $T_M^{\vee}$.

Under this Morita-equivalence the semi-simple $\hat{A}_{\mathfrak{m}}$-representation $M = S_1 \oplus \hdots \oplus S_n$ corresponds to the sum of the vertex-simples $\C e_1 \oplus \hdots \oplus \C e_n$, with the simple $S_i$ corresponding to the vertex-simple $\C e_i$ (the $e_i$ are the vertex-idempotents in the path algebra). Hence, also by Morita-equivalence we have an isomorphism
\[
Ext^1_{\hat{A}_{\mathfrak{m}}}(S_i,S_j) \simeq Ext^1_{\widehat{\C T_M^{\vee}}}(\C e_i,\C e_j) \]
Finally, because all ext-information is preserved under completions, and, because we know from representation-theory that the dimension of the ext-space between two vertex-simples for any quiver ${Q}$, 
$dim_{\C}~Ext^1_{\C Q}(\C e_i,\C e_j)$ is equal to the number of arrows starting in vertex $v_i$ and ending in vertex $v_j$, we are done!

Clearly, computing all $Ext^1_A(S,S')$ can still be a laborious task. However, it was proved in \cite{LBqurves} that all these dimensions follow often from a finite set of calculations when $A$ is a smooth algebra. The component semigroup $\wis{comp}(A)$ is the set of all connected components of the schemes $\wis{rep}_n~A$, for all $n \in \N$, with addition induced by the direct sum of finite dimensional representations.

The {\em one quiver} of $A$, $\wis{one}(A)$ is a full subquiver of $Q_A$ with one simple representant for every component which is a generator of $\wis{comp}(A)$ (note that such generators are determined by the fact that the component consists entirely of simples). Now, if $S$ and $T$ are two finite dimensional $A$-representations belonging to the connected components $\alpha$ and $\beta$ in $\wis{comp}(A)$ then we can write for certain $a_i,b_i \in \N$
\[
\alpha = \sum a_i g_i \quad \text{and} \quad \beta = \sum b_i g_i \]
with the $g_i$ the generator components. Then, $\epsilon=(a_i)_i$ and $\eta = (b_i)_i$ are dimension vectors for the one quiver. The main result of \cite{LBqurves} asserts now that
\[
dim_{\C}~Ext^1_A(S,T) = - \chi_{\wis{one}(A)}(\epsilon,\eta) \]
so that all ext-dimensions, and hence all $\mathfrak{m}$-adic completions of $A$ can be deduced from knowledge of the one quiver.

\subsection{}  \label{complexmodular} We will now make all these calculations explicit in the case of prime interest to us, which is the modular group $\Gamma = PSL_2(\Z)$, that is, we will describe the dual coalgebra $(\C \Gamma)^o$, at least in principle. Because $\Gamma \simeq C_2 \ast C_3$ we have that the group-algebra is the free algebra product of two semi-simple group algebras
\[
\C \Gamma \simeq \C C_2 \ast \C C_3 \]
and as such is a smooth algebra. In fact, a far more general result holds : whenever $G$ is a {\em virtually free group} (that is $G$ contains a free subgroup of finite index), then the group algebra $\C G$ is smooth by \cite{LBqurves}.

If $V$ is an $n$-dimensional $\Gamma$ representation, we can decompose it into eigenspaces for the action of $C_2 = \langle u \rangle$ and $C_3 = \langle v \rangle$ (let $\rho$ denote a primitive third root of unity) :
\[
V_+ \oplus V_- = V_1 \oplus V_2 = V \downarrow_{C_2} = V = V \downarrow_{C_3} = W_1 \oplus W_2 \oplus W_3 = W_1 \oplus W_{\rho} \oplus W_{\rho^2}  \]
If the dimension of $V_i$ is $a_i$ and that of $W_j$ is $b_j$, we say that $V$ is a $\Gamma$-representation of {\em dimension vector} $\alpha = (a_1,a_2;b_1,b_2,b_3)$. Choosing a basis $B_1$ of $V$ wrt. the decomposition $V_1 \oplus V_2$ and a basis $B_2$ wrt. $W_1 \oplus W_2 \oplus W_3$, we can view the basechange matrix $B_1 \rTo B_2$ as an $\alpha$-dimensional representation $V_Q$ of the quiver ${Q} = {Q}_{\Gamma}$
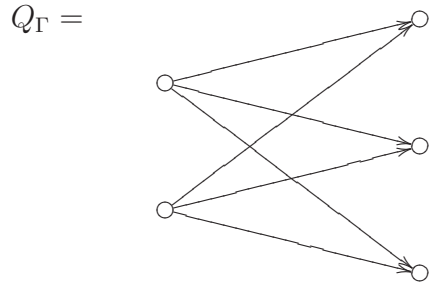
\[
{Q}_{\Gamma} = \qquad \xymatrix@=.6cm{
& & & & \vtx{} \\
\vtx{} \ar[rrrru] \ar[rrrrd] \ar[rrrrddd] & & & & \\
& & & & \vtx{} \\
\vtx{} \ar[rrrru] \ar[rrrruuu] \ar[rrrrd] & & & & \\
& & & &  \vtx{} }
\]
For a general quiver ${Q}$ on $k$ vertices, a weight $\theta \in \Z^k$ acts on the dimension vectors via the usual (Euclidian) scalar inproduct. A ${Q}$-representation of dimension vector $\alpha \in \N^k$ is said to be $\theta$-{\em stable} if and only if $\theta.\alpha=0$ and for every proper non-zero subrepresentation $W \subset V$ of dimension vector $\beta < \alpha$ we have that $\theta.\beta > 0$.

Bruce Westbury \cite{Westbury} has shown that $V$ is an irreducible $\Gamma$-representation if and only if $V_Q$ is a {\em $\theta$-stable} $Q$-representation where $\theta = (-1,-1;1,1,1)$ and that the two notions of isomorphism coincide. The {\em Euler-form} $\chi_{{Q}}$ of the quiver $Q$ is the bilinear form on $\mathbb{Z}^{\oplus 5}$ determined by the matrix
\[
\chi_{{Q}} = \begin{bmatrix} 1 & 0 & -1 & -1 & -1 \\ 0 & 1 & -1 & -1 & -1 \\ 0 & 0 & 1 & 0 & 0 \\ 0 & 0 & 0 & 1 & 0 \\ 0 & 0 & 0 & 0 & 1 \end{bmatrix} \]
Westbury also showed that if there exists a $\theta$-stable $\alpha$-dimensional $Q$-representation, then there is an $1 - \chi_{{Q}}(\alpha,\alpha)$ dimensional family of isomorphism classes of such representations (and a Zariski open subset of them will correspond to isomorphism classes of irreducible $\Gamma$-representations).

We will describe the one quiver $\wis{one}(\C \Gamma)$. By the above it follows that both the component semigroup $\wis{comp}~\C \Gamma$ and the semigroup of $\Z^5$ generated by all $\theta$-stable ${Q}$-representations are generated by the following six connected components, belonging to the dimension vectors
\[
g_{ij} = (\delta_{1i},\delta_{2i},\delta_{3i};\delta_{1j},\delta_{2j}) \]
and if we order and relabel these generators as
\[
a=g_{11}, b = g_{22}, c=g_{31}, d=g_{12}, e=g_{21}, f=g_{32} \]
we can compute from the Euler-form of ${Q}$ that the one-quiver of the modular group algebra is the following hexagonal graph
\[
\wis{one}(\C \Gamma) = \qquad 
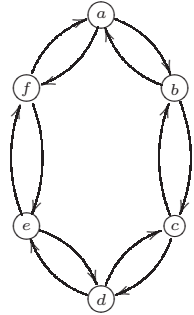
\xymatrix@=.6cm{
& \vtx{a} \ar@/^/[rd] \ar@/^/[ld] & \\
\vtx{f} \ar@/^/[ru] \ar@/^/[dd] & & \vtx{b} \ar@/^/[lu] \ar@/^/[dd] \\
& & \\
\vtx{e} \ar@/^/[uu] \ar@/^/[rd] & & \vtx{c} \ar@/^/[uu] \ar@/^/[ld] \\
& \vtx{d} \ar@/^/[lu] \ar@/^/[ru] &} \]
which is the origin of a lot of {\em hexagonal moonshine} in the representation theory of the modular group. In particular it follows from symmetry of the one quiver that the quiver $Q_{\C \Gamma}$ is also symmetric! 

\subsection{} Recall that an affine $\C$-algebra $A$ is said to be 2-Calabi-Yau if $gldim(A)=2$ and for any pair $S,T$ of finite dimensional $A$-representations, there exists a natural duality
\[
Ext^i_A(S,T) \simeq (Ext^{2-i}_A(T,S))^* \]
satisfying an additional sign condition. Raf Bocklandt \cite{Bocklandt} succeeded in extending the results on smooth algebras recalled before to the setting of 2-Calabi-Yau algebras. From the duality condition it is immediate that the quiver $Q_A$ is symmetric, that is, for every arrow $S \rTo^a T$ there is a paired arrow in the other direction $T \rTo^{a^*} S$. Bocklandt's result asserts that the $\mathfrak{m}$-adic completion $\hat{A}_{\mathfrak{m}}$ with $\mathfrak{m} = Ker(\pi_M)$ is Morita equivalent to the completion of the path algebra  of the (dual) quiver $Q_A | M$ modulo the preprojective relation
\[
\sum_a [a,a^*] = 0 \]
Further, he extends the idea of the one quiver to the 2-Calabi-Yau setting, allowing to compute the quiver $Q_A$ often from a finite set of calculations, using earlier results due to Crawley-Boevey \cite{CB}.

The group algebra $\C \pi_1(C)$ of the fundamental group of a genus $g$ Riemann surface is 2-Calabi-Yau by a result of Maxim Kontsevich. In \cite[\S 7.1]{Bocklandt} it is shown that the one-quiver of $\C \pi_1(C)$ consists of one vertex, corresponding to any one-dimensional simple representation, and $2g$ loops. From this and the results by Crawley-Boevey it follows that when $M=S_1^{\oplus e_1} \oplus \hdots \oplus S_k^{\oplus e_k}$ is a semi-simple $\C \pi_1(C)$-representation wit the simple factor $S_i$ having dimension $n_i$, then $\C Q_{\C \pi_1(C)} | M$ consists of $k$ vertices (corresponding to the distinct simple components $S_i$), such that the $i$-th vertex has exactly $2(g-1)n_i^2+2$ loops and there are exactly $2n_in_j(g-1)$ directed arrows from vertex $i$ to vertex $j$. 

This information allows us then to compute all $\mathfrak{m}$-adic completions of $\C \pi_1(C)$ as Morita equivalent to the completion of the path algebra of this quiver modulo the preprojective relation.

\subsection{}  In \ref{complexmodular} we described the structure of the path-coalgebra $\C Q_{\C \Gamma}$ which is Morita-Takeuchi equivalent to the dual complex coalgebra $(\C \Gamma)^o$. Describing the integral dual coalgebra $(\Z \Gamma)^o$ is a lot more complicated and will involve a good deal of knowledge of the integral (and modular) representation theory of the modular group.

Observe that the calculations in \ref{complexmodular} are valid for every algebraically closed field, so we might as well describe the coalgebra $(\overline{\Q} \Gamma)^o$ and study the action of the absolute Galois group $Gal(\overline{\Q}/\Q)$ on  it, giving us an handle on the rational dual coalgebra 
\[
(\Q \Gamma)^o = ((\overline{\Q} \Gamma)^o)^{Gal(\overline{\Q}/\Q)} \]
which brings us closer to $(\Z \Gamma)^o$. But, as in the case of the multiplicative group in the previous section, we do not require the full structure of this dual coalgebra but rather the image of the F-un data in it. As observed before, the $\mathbb{F}_1$-representation theory of $\Gamma$ is equivalent to the study of all finite dimensional transitive permutation representations of $\Gamma$ and hence to conjugacy classes of finite index subgroups of $\Gamma$.

We will recall the combinatorial description of those, due to R. Kulkarni in \cite{Kulkarni} in terms of generalized Farey symbols. Starting from this symbol we then describe how to associate a dessin, its monodromy group an finally to derive from it the modular content, that is the noncommutative gadget describing all $\Gamma$-representations deforming to the given permutation representation. In the next subsection we will give some interesting examples.

A {\em generalized Farey sequence} is an expression of the form
\[
\{ \infty=x_{-1},x_0,x_1,\hdots,x_n,x_{n+1}=\infty \} \]
where $x_0$ and $x_n$ are integers and some $x_i=0$. Moreover, all $x_i= \frac{a_i}{b_i}$ are rational numbers in reduced form and ordered such that
\[
| a_i b_{i+1} - b_i a_{i+1} | = 1 \qquad \text{for all $1 \leq i < n$} \]
The terminology is motivated by the fact that the classical {\em Farey sequence} $F(n)$, that is the ordered sequence of all rational numbers $0 \leq \frac{a}{b} \leq 1$ in reduced form with $b \leq n$, has this remarkable property.

A {\em Farey symbol} is a generalized Farey sequence $\{ \infty=x_{-1},x_0,x_1,\hdots,x_n,x_{n+1}=\infty \}$ such that for all $-1 \leq i \leq n$ we add one of the following symbols to two consecutive terms
\[
\xymatrix{x_i \ar@{-}[r]_{\bullet} & x_{i+1}} \quad \text{or} \quad \xymatrix{x_i \ar@{-}[r]_{\circ} & x_{i+1}}  \quad \text{or} \quad \xymatrix{x_i \ar@{-}[r]_{k} & x_{i+1}} \]
where each of the occurring  integers $k$ occur in pairs.

To connect Farey symbols with cofinite subgroups of the modular group $\Gamma$ we need to recall the {\em Dedekind tessellation} of the upper-half plane $\mathbb{H}$. Recall that the {\em extended modular group} $\Gamma^* = PGL_2(\mathbb{Z})$ acts on $\mathbb{H}$ via the natural action of $\Gamma$ on it together with the extra symmetry $z \mapsto - \overline{z}$. The Dedekind tessellation is the tessellation by fundamental domains for the action of $\Gamma^*$ on $\mathbb{H}$. It splits every fundamental domain for $\Gamma$ in two hyperbolic triangles, usually depicted as a black and a white one. Here is a depiction of the upper part of the Dedekind tessellation

\par \vskip 4mm

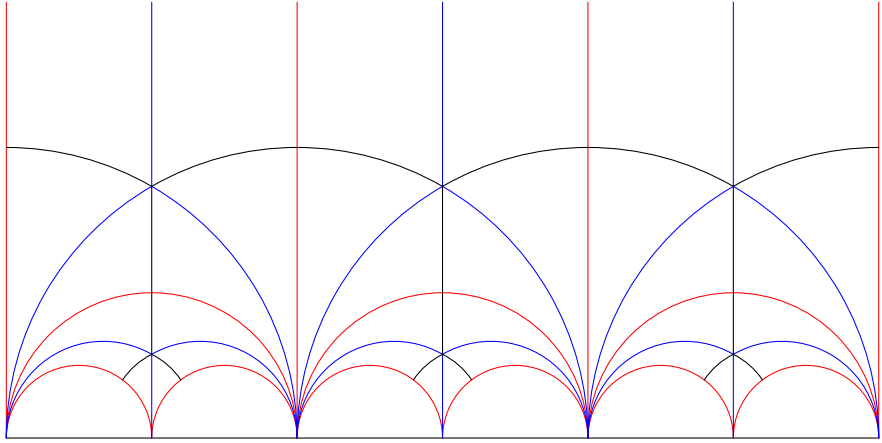
\begin{center}

\begin{mfpic}[110]{-1}{2}{0}{1.5}
\lines{(-1,0),(2,0)}

\drawcolor{red}
\lines{(0,0),(0,1.5)}
\lines{(-1,0),(-1,1.5)}
\lines{(1,0),(1,1.5)}
\lines{(2,0),(2,1.5)}
\arc[s]{(0,0),(-1,0),180}
\arc[s]{(1,0),(0,0),180}
\arc[s]{(2,0),(1,0),180}
\arc[s]{(-.5,0),(-1,0),180}
\arc[s]{(0,0),(-.5,0),180}
\arc[s]{(.5,0),(0,0),180}
\arc[s]{(1,0),(.5,0),180}
\arc[s]{(1.5,0),(1,0),180}
\arc[s]{(2,0),(1.5,0),180}

\drawcolor{blue}
\lines{(-.5,0.866),(-.5,1.5)}
\lines{(.5,0.866),(.5,1.5)}
\lines{(1.5,0.866),(1.5,1.5)}
\arc[s]{(-.5,0.866),(-1,0),60}
\arc[s]{(.5,0.866),(0,0),60}
\arc[s]{(1.5,0.866),(1,0),60}
\arc[s]{(0,0),(-.5,.866),60}
\arc[s]{(1,0),(.5,.866),60}
\arc[s]{(2,0),(1.5,.866),60}
\arc[s]{(-.5,.288),(-1,0),120}
\arc[s]{(.5,.288),(0,0),120}
\arc[s]{(1.5,.288),(1,0),120}
\arc[s]{(0,0),(-.5,.288),120}
\arc[s]{(1,0),(.5,.288),120}
\arc[s]{(2,0),(1.5,.288),120}
\lines{(-.5,0),(-.5,.288)}
\lines{(.5,0),(.5,.288)}
\lines{(1.5,0),(1.5,.288)}

\drawcolor{black}
\lines{(-.5,.866),(-.5,.288)}
\lines{(.5,.866),(.5,.288)}
\lines{(1.5,.866),(1.5,.288)}
\arc[s]{(-.5,.866),(-1,1),30}
\arc[s]{(.5,.866),(0,1),30}
\arc[s]{(1.5,.866),(1,1),30}
\arc[s]{(0,1),(-.5,.866),30}
\arc[s]{(1,1),(.5,.866),30}
\arc[s]{(2,1),(1.5,.866),30}
\arc[s]{(1.6,.2),(3/2,0.288),30}
\arc[s]{(0.6,.2),(1/2,0.288),30}
\arc[s]{(-.4,.2),(-1/2,0.288),30}
\arc[s]{(-.5,0.288),(-.6,.2),30}
\arc[s]{(.5,0.288),(.4,.2),30}
\arc[s]{(1.5,0.288),(1.4,.2),30}

\end{mfpic}

\end{center}
\par \vskip 4mm

Here, every red edge is a $\Gamma$-translate of the edge $[i,\infty]$, a blue edge a $\Gamma$-translate of $[\rho,\infty]$ where $\rho$ is a primitive sixth root of unity and every black edge is a $\Gamma$-translate of the circular arc $[i,\rho]$. Observe that every hyperbolic triangle of this tessellation has one edge of all three colors. Moving counterclockwise along the border of a triangle we either have the ordering red-blue-black (in which case we call this triangle a {\em white} triangle) or blue-red-black (and then we call it a {\em black} triangle). Any pair of a white and black triangle make a fundamental domain for the action of $\Gamma$.

Observe that any hyperbolic geodesic connecting two consecutive terms of a generalized Farey sequence consists of two red edges (connected at an intersection with black edges. We call these intersection points {\em even points} (later in the theory of dessins they will be denoted by a $\bullet$). A point where three blue edges come together with three black edges will be called an {\em odd point} (later denoted by $\xymatrix{\vtx{}}$).

A generalized Farey sequence therefore determines a hyperbolic polygonal region of $\mathbb{H}$ bounded by the (red) full geodesics connecting consecutive terms. The extra information contained in a Farey symbol tell us how to identify sides of this polygon (as well as how to extend it slightly in case of $\bullet$-connections) as follows :
\begin{itemize}
\item{For $\xymatrix{x_i \ar@{-}[r]_{\circ} & x_{i+1}}$ the two red edges making up the geodesic connecting $x_i$ with $x_{i+1}$ are identified.}
\item{For $\xymatrix{x_i \ar@{-}[r]_{k} & x_{i+1}}$ (with paired $\xymatrix{x_j \ar@{-}[r]_{k} & x_{j+1}}$) these two full geodesics (each consisting of two red edges) are identified.}
\item{For $\xymatrix{x_i \ar@{-}[r]_{\bullet} & x_{i+1}}$ we extend the boundary of the polygon by adding the two triangles just outside the full geodesic and identify the two blue edges forming the adjusted boundary.}
\end{itemize}
In this way, we associate to a Farey symbol a compact surface. Next, we will construct a {\em cuboid tree diagram} out of it, that is, a tree embedded in $\mathbb{H}$ such that all internal vertices are $3$-valent. Take as the vertices all odd-points lying in the interior of the polygonal region together with together with all even (red) and odd (blue) points on the boundary. We connect these vertices with the black lines in the interior of the polygonal region and add an involution on the red leaf-vertices determined by the side-pairing information contained in the Farey-symbol.

Finally, we will also associate to it a {\em bipartite cuboid graph} (aka a 'dessin d'enfants'). Start with the cuboid tree diagram and divide all edges in two (that is, add also the even internal points connecting the two black edges making up an edge in the tree diagram) and connect two red leaf-vertices when they correspond to each other under the involution.

For example, consider the Farey symbol
\[
\xymatrix{\infty \ar@{-}[r]_{1} & 0 \ar@{-}[r]_{\bullet} & \frac{1}{3} \ar@{-}[r]_{\bullet} & \frac{1}{2} \ar@{-}[r]_{\bullet} & 1 \ar@{-}[r]_{1} & \infty} \]
The boundary of the polygonal region determined by the symbol is indicated by the slightly thicker red and blue edges. The vertices of the cuboid tree are the red, blue and black points and the edges are the slightly thicker black edges.

\vskip 4mm

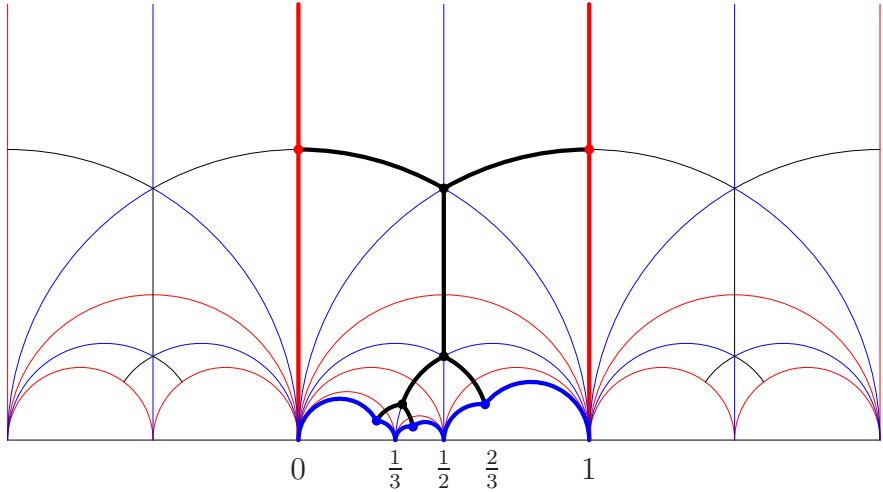
\begin{center}

\begin{mfpic}[110]{-1}{2}{0}{1.5}
\pen{.4}
\lines{(-1,0),(2,0)}

\drawcolor{red}
\lines{(-1,0),(-1,1.5)}
\pen{1.7}
\lines{(0,0),(0,1.5)}
\lines{(1,0),(1,1.5)}
\pen{.4}
\lines{(2,0),(2,1.5)}
\arc[s]{(0,0),(-1,0),180}
\arc[s]{(1,0),(0,0),180}
\arc[s]{(2,0),(1,0),180}
\arc[s]{(-.5,0),(-1,0),180}
\arc[s]{(0,0),(-.5,0),180}
\arc[s]{(.5,0),(0,0),180}
\arc[s]{(1,0),(.5,0),180}
\arc[s]{(1.5,0),(1,0),180}
\arc[s]{(2,0),(1.5,0),180}

\arc[s]{(1/3,0),(0,0),180}
\arc[s]{(1/2,0),(1/3,0),180}

\drawcolor{blue}
\lines{(-.5,0.866),(-.5,1.5)}
\lines{(.5,0.866),(.5,1.5)}
\lines{(1.5,0.866),(1.5,1.5)}
\arc[s]{(-.5,0.866),(-1,0),60}
\arc[s]{(.5,0.866),(0,0),60}
\arc[s]{(1.5,0.866),(1,0),60}
\arc[s]{(0,0),(-.5,.866),60}
\arc[s]{(1,0),(.5,.866),60}
\arc[s]{(2,0),(1.5,.866),60}
\arc[s]{(-.5,.288),(-1,0),120}
\arc[s]{(.5,.288),(0,0),120}
\arc[s]{(1.5,.288),(1,0),120}
\arc[s]{(0,0),(-.5,.288),120}
\arc[s]{(1,0),(.5,.288),120}
\arc[s]{(2,0),(1.5,.288),120}
\lines{(-.5,0),(-.5,.288)}
\lines{(.5,0),(.5,.288)}
\lines{(1.5,0),(1.5,.288)}

\arc[s]{(.357,.123),(1/3,0),18}

\pen{1.7}
\arc[s]{(0.268,0.0666),(0,0),152}
\arc[s]{(1/3,0),(0.268,0.0666),84}
\arc[s]{(0.394,0.046),(1/3,0),100}
\arc[s]{(1/2,0),(0.394,0.046),130}
\arc[s]{(0.642,.123),(1/2,0),95}
\arc[s]{(1,0),(0.642,.123),142}

\pen{.4}
\drawcolor{black}
\lines{(-.5,.866),(-.5,.288)}
\pen{1.7}
\lines{(.5,.866),(.5,.288)}
\pen{.4}
\lines{(1.5,.866),(1.5,.288)}
\arc[s]{(-.5,.866),(-1,1),30}
\pen{1.7}
\arc[s]{(.5,.866),(0,1),30}
\pen{.4}
\arc[s]{(1.5,.866),(1,1),30}

\arc[s]{(0,1),(-.5,.866),30}
\pen{1.7}
\arc[s]{(1,1),(.5,.866),30}
\pen{.4}
\arc[s]{(2,1),(1.5,.866),30}
\arc[s]{(1.6,.2),(3/2,0.288),30}
\arc[s]{(0.6,.2),(1/2,0.288),30}
\arc[s]{(-.4,.2),(-1/2,0.288),30}
\arc[s]{(-.5,0.288),(-.6,.2),30}
\arc[s]{(.5,0.288),(.4,.2),30}
\arc[s]{(1.5,0.288),(1.4,.2),30}

\pen{1.7}
\arc[s]{(.5,.288),(0.357,.123),45}
\arc[s]{(0.357,0.123),(0.269,0.0666),55}
\arc[s]{(0.394,0.046),(0.357,0.123),30}
\arc[s]{(.642,.123),(.5,.288),47}

\drawcolor{blue}
\pen{2}
\circle{(0.269,0.0666),.008}
\circle{(0.394,0.046),.008}
\circle{(0.642,.123),.008}

\drawcolor{red}
\circle{(0,1),.008}
\circle{(1,1),.008}

\drawcolor{black}
\circle{(0.357,0.123),.008}
\circle{(1/2,0.288),.008}
\circle{(1/2,0.866),.008}

\tlabel[cc](.5,-.1){$\frac{1}{2}$}
\tlabel[cc](.333,-.1){$\frac{1}{3}$}
\tlabel[cc](.666,-.1){$\frac{2}{3}$}
\tlabel[cc](0,-.1){$0$}
\tlabel[cc](1,-.1){$1$}

\end{mfpic}

\end{center}

\vskip 4mm

Because the two red leaf-vertices correspond to each other under the involution, the corresponding bipartite cuboid diagram (or modular dessin) is 
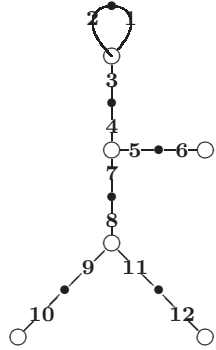
\[
\xymatrix@=1cm{
& \vtx{} \aar{1}{2}@(ur,ul) \aar{3}{4}[d] & \\
& \vtx{} \aar{5}{6}[r] \aar{7}{8}[d] & \vtx{} \\
& \vtx{} \aar{9}{10}[ld] \aar{11}{12}[rd] & \\
\vtx{} & & \vtx{} }
\]

Such a dessin encodes the data of a Belyi covering $C \rOnto \mathbb{P}^1_{\C}$ ramified only in the points $\{ 0,1,\infty \}$.
The inverse images of $0$ will be represented by a $\xymatrix{\vtx{}}$-vertex, those of $1$ by a $\bullet$-vertex. Of relevance for us are dessins which are {\em modular quilts} meaning that every $\bullet$-vertex is $2$-valent and every $\xymatrix{\vtx{}}$-vertex is $1$- or $3$-valent.

Given a modular dessin, denote each of the edges by a different number between $1$ and $d$ (the degree of $\pi$), then the {\em monodromy group} $G_{\pi}$ of $\pi$ is the subgroup of $S_d$ generated by the order three element $\sigma_0$ obtained by cycling round every $\xymatrix{\vtx{}}$
-vertex counterclockwise and the order two element $\sigma_1$ obtained by recording the two edges ending at every $\bullet$-vertex. This defines an exact sequence of groups
\[
1 \rTo G \rTo \Gamma \rTo G_{\pi} \rTo 1 \]
and the projective curve $C$ corresponding to the modular dessin can be identified with a compactification of $\mathbb{H} / G$ where $\mathbb{H}$ is the upper half-plane on which $G \subset \Gamma$ acts via M\"obius transformations.

The $d$-dimensional permutation representation $M =  \Gamma/G$ decomposes into irreducible representations for the monodromy group $G_{\pi}$, say
\[
M = X_1^{\oplus e_1} \oplus \hdots \oplus X_k^{\oplus e_k} \]
with every $X_i$ an irreducible $G_{\pi}$ and hence also irreducible $\Gamma$-representation. The {\em modular content} of the dessin, or of the permutation representation, is the quiver on $k$ vertices $Q_{\pi} = Q_{\C \Gamma} | M$ together with the dimension vector $\alpha_{\pi} = (e_1,\hdots,e_k)$ determined by the multiplicities of the simples in the permutation representation.

Roughly speaking, the modular content $(Q_{\pi},\alpha_{\pi})$ encodes how much the curve $C$, the dessin or the permutation representation 'sees' of the modular group. That is, the quotient variety $\wis{iss}_{\alpha_{\pi}} Q_{\pi} = \wis{rep}_{\alpha_{\pi}} Q_{\pi} / GL(\alpha_{\pi})$ classifies all semi-simple $d$-dimensional $\Gamma$-representations deforming to the permutation representation $M$. As such, it is a new noncommutative gadget associated to a classical object, the curve $C$. It would be interesting to know whether the modular content is a Galois invariant of the dessin, or more generally, what subsidiary information derived from it is a Galois invariant.

We now give an algorithm to compute the modular content, using the group-theory program GAP, starting from the modular quilt $D$.

\begin{enumerate}
\item{Determine the permutations $\sigma_0,\sigma_1 \in S_d$ described above, that is obtained by walking around the $\bullet$-vertices (for $\sigma_1$) and the $\xymatrix{\vtx{}}$-vertices (for $\sigma_0$) in $D$ and feed them to GAP as \texttt{s0,s1}.}
\item{Calculate the monodromy group $G_{\pi}$ via \texttt{G:=Group(s0,s1)} and determine its character table via \texttt{chars:=CharacterTable(G);)}}
\item{Determine the $G_{\pi}$-character of the permutation representation by calling \texttt{ConjugacyClasses(G)}. This returns a list of $S_d$-permutations representing the conjugacy classes of $G_{\pi}$. To determine the character-value we only need to count the numbers missing in the cycle decomposition of the permutation. Let $\chi$ be the obtained character which is the list \texttt{chi}.}
\item{Determine the irreducible components of $\chi$ and their multiplicities via
\texttt{MatScalarProducts(chars,Irr(chars),[chi]);}. The non-zero entries form the dimension vector $\alpha_{\pi}$ and they determine the simple factors $X_1,\hdots,X_k$.}
\item{Determine the conjugacy classes of $\sigma_0$ and $\sigma_1$. For example,  the number of the conjugacy class in the character table is found by \texttt{FusionConjugacyClasses(Group(s0),G);}. Alternatively, one can use \texttt{IsConjugate(G,s0,s);} for \text{s} a suitable element representant obtained via \texttt{ConjugacyClasses(G);}. Assume $\sigma_0$ (resp. $\sigma_1$) belongs to the $a$-th (resp. $b$-th) conjugacy class.}
\item{From the character values of $X_i$ in the $a$-th and $b$-th column of \texttt{Display(chars);} one deduces the dimension vector $\alpha_i=(a_1(i),a_2(i);b_1(i),b_2(i),b_3(i))$ of the ${Q}_{\Gamma}$-representation corresponding to $X_i$.}
\item{Finally, the number of arrows (and loops) in the quiver ${Q}_{\pi}$ between the vertices corresponding to $X_i$ and $X_j$ is given by $\delta_{ij}-\chi_{{Q}_{\Gamma}}(\alpha_i,\alpha_j)$.}
\end{enumerate}

\subsection{}  As the modular content encodes all possible $\Gamma$-representation deformations of the permutation representation, it is often a huge object which makes it difficult to extract interesting deformations from it. Sometimes though, a true gem reveals itself. In the previous subsection we used  the generalized Farey-symbol
\[
\xymatrix{\infty \ar@{-}[r]_{1} & 0 \ar@{-}[r]_{\bullet} & \frac{1}{3} \ar@{-}[r]_{\bullet} & \frac{1}{2} \ar@{-}[r]_{\bullet} & 1 \ar@{-}[r]_{1} & \infty} \]
 Note that it consists of half of the Farey-sequence $F(3)$ (those $\leq \frac{1}{2}$). Generalizing this construction for all classical Farey sequences leads to an intriguing class of examples. The $n$-th {\em Iguanodon Farey-symbol} is the Farey symbol
\[
\xymatrix{\infty \ar@{-}[r]_{1} & 0 \ar@{-}[r]_{\bullet} & \frac{1}{n} \ar@{-}[r]_{\bullet} & \hdots & \frac{1}{2} \ar@{-}[l]^{\bullet} \ar@{-}[r]_{\bullet} & 1 \ar@{-}[r]_{1} & \infty} \]
where the rational numbers occurring are precisely those Farey numbers in $F(n)$ smaller or equal to $\frac{1}{2}$.

The terminology is explained by depicting the first few bipartite cuboid diagrams associated to Farey sequences
{\tiny
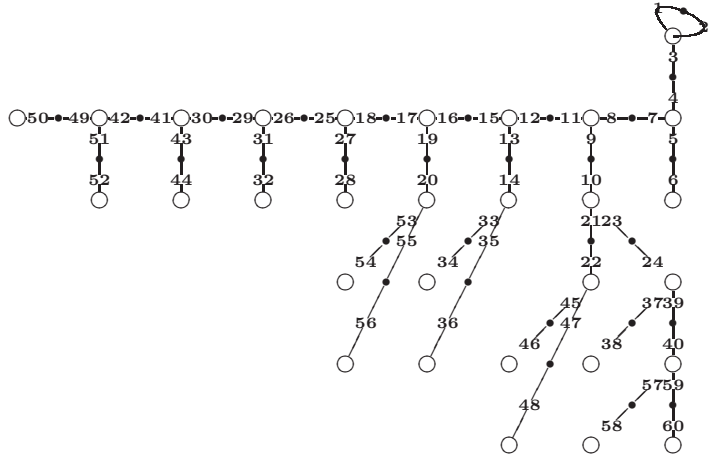
\[
\xymatrix{
& & & & & & & & \vtx{} \aar{3}{4}[d] \aar{1}{2}@(ul,r) \\
\vtx{} \aar{50}{49}[r] & \vtx{} \aar{51}{52}[d] \aar{42}{41}[r] & \vtx{} \aar{43}{44}[d] \aar{30}{29}[r] & \vtx{} \aar{26}{25}[r] \aar{31}{32}[d] & \vtx{} \aar{18}{17}[r] \aar{27}{28}[d] & \vtx{} \aar{16}{15}[r] \aar{19}{20}[d] & \vtx{} \aar{12}{11}[r] \aar{13}{14}[d] & \vtx{} \aar{8}{7}[r] \aar{9}{10}[d] & \vtx{} \aar{5}{6}[d] \\
& \vtx{} & \vtx{} & \vtx{} & \vtx{} & \vtx{} & \vtx{} & \vtx{} & \vtx{} \\
& & & & \vtx{} \aar{54}{53}[ru] & \vtx{} \aar{34}{33}[ru] & & \vtx{} \aar{22}{21}[u] & \vtx{} \aar{24}{23}[lu] \\
& & & & \vtx{} \aar{56}{55}[ruu] & \vtx{} \aar{36}{35}[ruu] & \vtx{} \aar{46}{45}[ru] & \vtx{} \aar{38}{37}[ru] & \vtx{} \aar{40}{39}[u] \\
& & & & & & \vtx{} \aar{48}{47}[ruu] & \vtx{} \aar{58}{57}[ru] & \vtx{} \aar{60}{59}[u]}
\]
}
Here, the diagram corresponding to Farey sequence $F(n)$ is the full subfigure on the first $m(n)$ (half)edges
\[
\begin{array}{c|cccccccc}
n & 2 & 3 & 4 & 5 & 6 & 7 & 8 & 9 \\
\hline 
m(n) & 8 & 12 & 16 & 24 & 28 & 40 & 48 & 60
\end{array}
\]
The monodromy groups corresponding to the $n$-th Iguanodon symbol are
\[
\begin{array}{c|cccccccc}
n & 2 & 3 & 4 & 5 & 6 & 7 & 8 & 9 \\
\hline
& L_2(7) & M_{12} & A_{16} & M_{24} & A_{28} & A_{40} & A_{48} & A_{60} \\
\hline \\ \hline
n & 10 & 11 & 12 & 13 & 14 & 15 & 16 & 17 \\
\hline
& A_{68} & A_{88}& A_{96} & A_{120} & A_{132} & A_{148} & A_{164} & A_{196}
\end{array}
\]
This can be verified by hand (and GAP) using the above picture for $n \leq 9$ and by using the SAGE-package \texttt{kfarey.sage} for higher $n$. It is plausible that the monodromy groups of the Iguanodon symbols are all simple groups and it is quite remarkable that the Mathieu groups $M_{12}$ and $M_{24}$ appear in this sequence of alternating groups.

Now, let us compute the modular content of these permutation representations. The action of the monodromy group is clearly 2-transitive implying that as a $\C G_{\pi}$-representation, the permutation representation splits into two irreducibles, one of which being clearly the trivial representation. Note also that the character of the generator of order $2$ is equal to zero as there are no $\bullet$-end points. Further, $\circ$-endpoints appear in pairs and add another 4 half-edges, that is 4 dimensions, to the permutation space. By induction we see that the dimension of the permutation representation is always of the form $4n$ with $\chi(\sigma_1)=0$ and $\chi(\sigma_0)=n$. 

By the argument recalled in \ref{complexmodular} it follows that the dimension vector of the $Q_{\Gamma}$-quiver representation corresponding to the permutation representation is
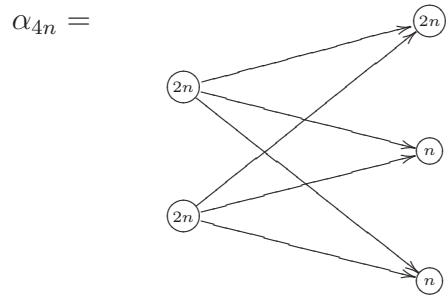
\[
\alpha_{4n}= \qquad \xymatrix@=.5cm{
& & & & \vtx{2n} \\
\vtx{2n} \ar[rrrru] \ar[rrrrd] \ar[rrrrddd] & & & & \\
& & & & \vtx{n} \\
\vtx{2n} \ar[rrrru] \ar[rrrruuu] \ar[rrrrd] & & & & \\
& & & &  \vtx{n} }
\]
By 2-transitivity the dimension vectors of the two simple components $S$ and $T$ are 
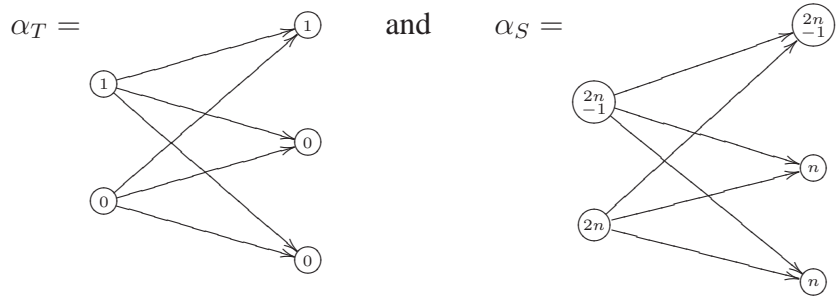
\[
\alpha_T=\xymatrix@=.4cm{
& & & & \vtx{1} \\
\vtx{1} \ar[rrrru] \ar[rrrrd] \ar[rrrrddd] & & & & \\
& & & & \vtx{0} \\
\vtx{0} \ar[rrrru] \ar[rrrruuu] \ar[rrrrd] & & & & \\
& & & &  \vtx{0} } \qquad \text{and} \qquad
\alpha_S= \xymatrix@=.4cm{
& & & & \vtx{\overset{2n}{-1}} \\
\vtx{\overset{2n}{-1}} \ar[rrrru] \ar[rrrrd] \ar[rrrrddd] & & & & \\
& & & & \vtx{n} \\
\vtx{2n} \ar[rrrru] \ar[rrrruuu] \ar[rrrrd] & & & & \\
& & & &  \vtx{n} }
\]
But then, by the algorithm we have that the modular content $(Q_{\pi},\alpha_{\pi})$ of the permutation representation can be depicted as

\vskip 4mm

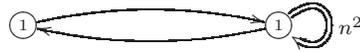
\[
\xymatrix{\vtx{1} \ar@/^/[rrr] & & & \vtx{1} \ar@/^/[lll] \ar@{=>}@(ur,dr)^{n^2}}
\]

\vskip 4mm

The $n^2$ loops in the vertex corresponding to the simple factor $S$ indicate that the moduli space of semi-stable $Q_{\Gamma}$-representations $M_{\theta}^{ss}(Q_{\Gamma},\alpha_S)$ is $n^2$-dimensional and as $S$ is a smooth point in it, there is an $n^2$-dimensional family of simple $\Gamma$-representations in the neighborhood of $S$. More interesting is the fact that there is just one arrow in each direction between the two vertices.

This implies that the permutation representation is a smooth point in the moduli space of semi-simple $\Gamma$-representations, a rare fact for higher dimensional decomposable representations (see the paper \cite{LBBocklandtSymens} for more details on singularities of quiver-representations). Further, this implies that there is a unique (!) curve of simple $4n$-dimensional $\Gamma$-representations degenerating to the given permutation representation! Certainly in the case of the sporadic Mathieu groups it would be interesting to study these curves (and their closures in the moduli space $M^{ss}_{\theta}(Q_{\Gamma},\alpha_{4n})$) in more detail.


\begin{thebibliography}{10}
 
 \bibitem{Bocklandt} Raf Bocklandt, {\it Noncommutative tangent cones and Calabi-Yau algebras}, arXiv:0711.0179 (2007)
 
 \bibitem{LBBocklandt} Raf Bocklandt and Lieven Le Bruyn, {\it Necklace Lie algebras and noncommutative symplectic geometry}, Math. Z. 240 (2002) 141-167, arXiv:math/0010030
 
 \bibitem{LBBocklandtSymens} Raf Bocklandt, Lieven Le Bruyn and Geert Van de Weyer, {\it Smooth order singularities}, J. Alg. Appl. 2 (2003) 365-395, arXiv:math/0207250
 
 \bibitem{CC} Alain Connes and Katia Consani, {\it On the notion of geometry over $\mathbb{F}_1$}, arXiv:0809.2926

\bibitem{CB} Bill Crawley-Boevey, {\it Geometry of the moment map for representations of quivers}, Compositio Math. 126 (2001) 257-293
 
 \bibitem{CuntzQuillen}
 Joachim Cuntz, Daniel  Quillen, {\it Algebra extensions and nonsingularity}, Journal of AMS, 
v.8, no. 2 (1995) 251Ð289 

\bibitem{Habiro} Kazuo Habiro, {\it Cyclotomic completions of polynomial rings}, arXiv:0209324



\bibitem{KontSoib}
Maxim Kontsevich and Yan Soibelman, {\it Notes on $A_{\infty}$-algebras, $A_{\infty}$-categories and non-commutative geometry I}, arXiv:math.RA/0606241 (2006)

\bibitem{Kulkarni}
Ravi S. Kulkarni, {\it An arithmetic-geometric method in the study of the subgroups of the modular group}, Amer. J. Math. 113 (1991) 1053-1133

\bibitem{LBqurves}
Lieven Le Bruyn, {\it Qurves and quivers}, arXiv:math.RA/0406618 (2004),
Journal of Algebra 290 (2005) 447-472

\bibitem{LBdualcoalgebra}
Lieven Le Bruyn, {\it Noncommutative geometry and dual coalgebras}, arXiv:0805.2377v1 (2008)


\bibitem{Manin} Yuri I. Manin, {\it Cyclotomy and analytic geometry over $\mathbb{F}_1$}, arXiv:0809.1564 (2008)

\bibitem{Marcolli} Matilde Marcolli, {\it Cyclotomy and endomotives}, arXiv:0901.3167 (2009)

\bibitem{Soule} Christophe SoulŽ, {\it Let vari\'et\'es sur le corps \`a un \'el\'ement}, Moscow Math. J. 4 (2004) 217-244


\bibitem{Sweedler}
Moss E. Sweedler, {\it Hopf Algebras}, monograph, W.A. Benjamin (New York) (1969)

\bibitem{Westbury}
Bruce Westbury, {\it On the character varieties of the modular group}, preprint Nottingham (1995)

\bibitem{Zagier} Don Zagier, {\it Vassiliev invariants and a strange identity related to the Dedekind eta-function}, available as MPI-preprint 

\end{thebibliography}
\end{document}